\documentclass[12pt,a4paper]{article}
\usepackage[dvipdfmx]{graphicx}
\usepackage{latexsym}
\usepackage{amssymb}
\usepackage{amsmath}
\usepackage{amscd}
\usepackage{amsfonts}
\usepackage{mathrsfs}
\usepackage{enumerate}
\usepackage{bm}
\usepackage[all]{xy}

\newtheorem{Th}{Theorem}[section]

\newtheorem{Pro}[Th]{Proposition}
\newtheorem{Rem}[Th]{Remark}
\newcommand{\demo}{\par\noindent{\it Proof. \/}\ }
\newcommand{\enD}{\hfill $\Box$\vspace{3truemm} \par}
\newcommand{\bx}{\mbox{\boldmath $x$}}
\newcommand{\bX}{\mbox{\boldmath $X$}}

\newcommand{\be}{\mbox{\boldmath $e$}}
\newcommand{\ba}{\mbox{\boldmath $a$}}
\newcommand{\bb}{\mbox{\boldmath $b$}}
\newcommand{\bv}{\mbox{\boldmath $v$}}
\newcommand{\by}{\mbox{\boldmath $y$}}

\newcommand{\bw}{\mbox{\boldmath $w$}}

\newcommand{\bxi}{\mbox{\boldmath $\xi$}}

\newcommand{\bsigma}{\mbox{\boldmath $\sigma$}}
\newcommand{\bgamma}{\mbox{\boldmath $\gamma$}}
\newcommand{\bt}{\mbox{\boldmath $t$}}
\newcommand{\bn}{\mbox{\boldmath $n$}}

\newcommand{\bc}{\mbox{\boldmath $c$}}

\newcommand{\sbgamma}{\mbox{\scriptsize \boldmath$\gamma$}}

\newcommand{\kn}{\kappa_n(s)}
\newcommand{\kg}{\kappa_g(s)}
\newcommand{\tg}{\tau_g(s)}
\newcommand{\dtr}{\delta^{T}_{r}}
\newcommand{\dsr}{\delta^{S}_{r}}
\newcommand{\dso}{\delta^{S}_{o}}
\newcommand{\dlo}{\delta^{L}_{o}}
\newcommand{\dlr}{\delta^{L}_{r}}
\newcommand{\ts}{\boldsymbol{t}}

\newcommand{\ns}{\boldsymbol{n}_\gamma}
\newcommand{\kns}{\kappa_n}
\newcommand{\kgs}{\kappa_g}
\newcommand{\tgs}{\tau_g}

\newcommand{\R}{{\mathbb R}}
\newcommand{\lon}{\longrightarrow}
\begin{document}
\title{Lorentzian Darboux images of curves on spacelike surfaces in Lorentz-Minkowski $3$-space }
\author{Noriaki ITO and Shyuichi IZUMIYA}

\date{\today}
\maketitle
\begin{abstract}
For a regular curve on a spacelike surface in Lorentz-Minkowski $3$-space, we have a moving frame along the curve which is called a Lorentzian Darboux frame.
We introduce five special vector fields along the curve associated to the Lorentzian Darboux frame and investigate their singularities.
\end{abstract}
\renewcommand{\thefootnote}{\fnsymbol{footnote}}
\footnote[0]{2010 Mathematics Subject classification. Primary 53A04;
Secondary 58Kxx} \footnote[0]{Key Words and Phrases. Lorentzian Darboux frame, Lorentzian Darboux vector field, pseudo-spherical image, singularities} 

\section{Introduction}
In this paper we consider a curve on a spacelike surface in the Lorentz-Minkowski $3$-space and some special vector fields along the curve.
The study of geometry of the Lorentz-Minkowski space is of interest in the special relativity theory. 
From the view point of mathematics, the interesting problem is how geometric properties of the Lorentz-Minkowski space is different from
those of the Euclidean space.
In the Euclidean $3$-space, the notion of Darboux frames along  curves on surfaces is well-known.
In \cite{HII} spherical duals (cf. \cite{Arnold2, Porteous}) of basis of the Darboux frame along a curve are introduced, which are called
Darboux vector fields along the curve.
There are three Darboux vector fields along the curve. Singularities and geometric properties of these three Darboux vectors were
investigated in \cite{HII}.
\par
On the other hand, there also exists a Lorentzian version of Darboux frames along curves on spacelike surfaces \cite{Sato}.
We consider (pseudo-spherical) Legendrian duals (cf. \cite{CIZ09,Izu2}) of basis of the Lorentzian Darboux frame along a curve, which are called Lorentzian Darboux vectors along the curve.
Since there are three kinds of pseudo-spheres in Lorentz-Minkowski space, we have eight Lorentzian Darboux vectors along the curve.
There are three Legendrian duals of the unit tangent vector along the curve, which were essentially investigated in \cite{Sato}. Those vector fields are
three of the Lorentzian Darboux vector fields along the curve.
Therefore, we consider remaining five Lorentzian Darboux vectors along the curve here. 
We investigate the singularities of the pseudo-spherical image of Lorentzian Darboux vectors.
As a consequence, we obtain five new Lorentzian invariants which characterize the singularities of these Lorentzian Darboux vectors.
We also investigate the geometric meanings of these invariants.

\section{Basic concepts}d
In this section we prepare some definitions and basic facts which we will use in this paper. For basic concepts and details of properties, see \cite{Oneil, Sato}.
Let $\R ^3$ be a three-dimensional vector space. 
For any $\bx=(x_0,x_1,x_2),\by=(y_0,y_1,y_2) \in \R^3$, the pseudo-scalar product of $\bx$ and $\by$ is defined by $\langle \bx,\by \rangle =-x_0y_0+x_1y_1+x_2y_2$. We call $(\R^3 , \langle,\rangle)$ the {\it Lorentz-Minkowski $3$-space\/}. We write $\R_1^3$ instead of $(\R^3 , \langle,\rangle).$ We say that a non-zero vector $\bx \in \R_1^3$ is {\it spacelike, lightlike\/} or {\it timelike} if $\langle \bx,\bx\rangle >0$ , $\langle \bx,\bx \rangle=0$ or $\langle \bx,\bx \rangle<0$, respectively. The norm of the vector $\bx \in \R_1^3$ is defined by $\parallel\bx\parallel =\sqrt{\lvert\langle \bx,\bx \rangle \rvert}$.  For a non-zero vector $\bv \in \R_1^3$ and a real number $c\in \R,$ we define a {\it plane} with a {\it pseudo-normal} $\bv$ by
$$
P(\bv,c)=\{\bx\in \R_1^3 \ | \ \langle \bx,\bv \rangle=c\ \}.
$$
We call $P(\bv,c)$ a {\it spacelike plane}, a {\it timelike plane} or a {\it lightlike \it plane} if $\bv$  is timelike, spacelike or lightlike, respectively. 
We introduce three pseudo-spheres in $\R^3_1$ as follows:
We define the {\it hyperbolic plane} by
$$
H^2(-1)=\{\bx\in \R^{3}_1 \ | \  \langle \bx ,\bx\rangle =-1
\},
$$
{\it de Sitter $2$-space} by
$$
S^2_1=\{\bx\in \R^{3}_1 \ | \ \langle \bx ,\bx\rangle =1\ \}
$$
and the {\it {\rm (}open{\rm )} lightcone} by
\[ LC^* = \{\bx \in \R^{3}_1\backslash \{ \boldsymbol{0} \} \ | \ \langle \bx,\bx \rangle = 0 \  \}.
\]
We also define the following curves on the pseudo-spheres with constant curvatures:
A curve defined by the intersection of $H^2(-1)$ with a plane is called a {\it hyperbolic line} (respectively, a {\it horocycle}) if 
the plane is a timelike plane through the origin (respectively, a lightlike plane).
We also say that a curve on the de Sitter $2$-space $S^2_1$ is a {\it geodesic pseudo-circle} (respectively, a {\it geodesic hyperbola})
if it is defined by the intersection of $S^2_1$ with a spacelike (respectively, a timelike) plane through the origin.
Moreover, a curve on $S^2_1$ is said to be a {\it de Sitter horocycle} if it is defined by the intersection of $S^2_1$ with a lightlike plane away from the origin.
Here we define 
\[
\ba\wedge \bb=
\begin{vmatrix}
-\be_0 & \be_1 & \be_2 \cr
a_0 & a_1 & a_2 \cr
b_0 & b_1 & b_2 
\end{vmatrix},
\]
where $\ba =(a_0,a_1,a_2), \bb=(b_0,b_1,b_2)$ and $\{\be _0,\be_1,\be_2\}$ is the canonical basis of $\R^3.$
\par
\par

We now prepare some basic facts of curves on a spacelike surface. We consider a spacelike embedding $\bX:U\lon
\R^{3}_1$ from an open subset $U\subset \R^2.$ We write
$M=\bX(U)$ and identify $M$ and $U$ through the embedding $\bX.$ 
Here, we say that $\bX$ is a {\it spacelike embedding}
if the tangent space $T_pM$ consists of spacelike vectors at any $p=\bX(u).$
Let $\bar{\bgamma}:I\lon U$ be a regular curve and we have a curve  $\bgamma:I\lon M\subset \R^3_1$ defined by $\bgamma(s)=\bX(\bar{\bgamma}(s)).$ We say that $\bgamma$ is a {\it curve on the spacelike surface M.} Since  $\bgamma$ is a spacelike curve, we can reparametrize it by the arc-length $s$.
So we have the spacelike unit tangent vector $\bt(s)=\bgamma'(s)$ of $\bgamma(s)$.
Since $\bX$ is a spacelike embedding, we have a timelike unit normal vector field $\bn$ along $M=\bX(U)$
defined by
\[
\bn(p)=\frac{\bX_{u_1}(u)\wedge \bX_{u_2}(u)}{\|\bX_{u_1}(u)\wedge \bX_{u_2}(u)\|},
\]
for $p=\bX(u)$.
We say that $\bn$ is {\it future directed}
if $\langle \bn ,\be _0\rangle <0.$
We choose the orientation of $M$ such that $\bn$ is future directed.
We define $\bn_{\sbgamma}(s)=\bn\circ\bgamma (s),$ so that we have a timelike unit normal vector
field $\bn_{\sbgamma}$ along $\bgamma .$ 
Therefore we can construct  a spacelike unit normal vetoer field $\bb(s) \in N_p(M)$ defined by
$\bb(s)=\bt(s)\wedge\bn_{\gamma}(s)$.
It follows that we have $\langle \bn_\gamma,\bn_\gamma\rangle =-1,\ \langle
\bn_\gamma,\bb\rangle =0,\ \langle \bb,\bb\rangle =1.$ Then we have a pseudo-orthonormal frame $\{\bt(s), \bn_\gamma(s), \bb(s)\}$ along $\bm{\gamma}$,
which is called the {\it Lorentzian Darboux frame} along $\bgamma$. By standard arguments, we have the following {\it Frenet-Serret type formulae}:
   $$\left\{
       \begin{array}{ll}
          \bt'(s)=\kappa_n(s)\bn_\gamma(s)+\kappa_g(s)\bb(s),  \\
          \bn_\gamma'(s)=\kappa_n(s)\bt(s)+\tau_g(s)\bb(s),\\
          \bb'(s)=-\kappa_g(s)\bt(s)+\tau_g(s)\bn_\gamma(s),\\
   \end{array}
     \right.$$   
where $\kappa_n(s)= -\langle \bt'(s),\bn_\gamma(s)\rangle$, $\kappa_g(s)= \langle \bt'(s),\bb(s)\rangle$ and $\tau_g(s)= -\langle \bb'(s),\bn_\gamma(s)\rangle.$
We have the geometric characterizations of $\bgamma$ by the invariants 
$\kappa_g, \kappa_n$ and $\tau_g.$
We say that $\bm{\gamma}$ is a {\it geodesic curve} if the curvature vector $\bm{t}'(s)$ has only a pseudo-normal 
component of the surface, an {\it asymptotic curve} if $\bm{t}'(s)$ has only a tangential component of the surface
and a {principal curve} if $\bm{n}_{\bm{\gamma}}'(s)$ is equal to the tangent direction of $\bm{\gamma}$, respectively. Then
 $$\bgamma \ \ \mbox{is} \ \ 
  \left\{ 
  \begin{array}{ll}
      \mbox{a geodesic curve if and only if}\ \kappa_g \equiv 0,\\
     \mbox{an asymptotic curve if and only if}\  \kappa_n \equiv 0,\\
      \mbox{a principal curve if and only if}\ \tau_g \equiv 0. \\  
   \end{array}
     \right.$$   
Then we define the following five {\it pseudo-spherical Lorentzian Darboux images\/} along $\bm{\gamma}$:
\begin{align*}
(A)& \ \overline{\bm{D}}^T_r:I\lon H^2(-1);\ \overline{\bm{D}}^T_r(s) = \frac{\tg\bt(s)-\kg\bn_{\bm{\gamma}}(s)}{\sqrt{\kg^2 - \tg^2}}\ \mbox{if}\ \kappa_g(s)^2>\tau_g(s)^2,\\           
(B)& \ \overline{\bm{D}}^S_r:I\lon S^2_1;\ \overline{\bm{D}}^S_r(s) = \frac{\tg\bt(s)-\kg\bn_{\bm{\gamma}}(s)}{\sqrt{\tg^2 - \kg^2}}\ \mbox{if}\ \tau_g(s)^2>\kappa_g(s)^2,\\           
(C)& \ \overline{\bm{D}}^L_r:I\lon LC^*;\ \overline{\bm{D}}^L_r(s) =  \frac{\tg\bt(s)-\kg\bn_{\bm{\gamma}}(s)}{\sqrt{\kg^2 - \tg^2}}  + \bb(s)\ \mbox{if}\ \kappa_g(s)^2>\tau_g(s)^2,\\ 
(D)& \ \overline{\bm{D}}^S_o:I\lon S^2_1;\ \overline{\bm{D}}^S_o(s) = \frac{\tg\bt(s)-\kn\bb(s)}{\sqrt{\kn^2 + \tg^2}}\ \mbox{if}\ (\kappa _n(s),\tau _g(s))\not= (0,0), \\          
(E)& \ \overline{\bm{D}}^L_o:I\lon LC^*;\ \overline{\bm{D}}^L_o(s) = \frac{\tg\bt(s)-\kn\bb(s)}{\sqrt{\kn^2 + \tg^2}}  + \bn_{\bm{\gamma}}(s)\ \mbox{if}\ (\kappa _n(s),\tau _g(s))\not= (0,0). 
\end{align*}
We call $(A)$ the {\it pseudo-spherical rectifying timelike Darboux image}, $(B)$ the {\it pseudo-spherical rectifying spacelike Darboux image}, $(C )$ the {\it pseudo-spherical rectifying lightlike Darboux image},
$(D)$ the {\it pseudo-spherical osculating spacelike Darboux image} and $(E)$ the {\it pseudo-spherical osculating lightlike Darboux image} along $\bgamma ,$
respectively.
We remark that we cannot define a pseudo-spherical osculating timelike Darboux image.
\begin{Rem}
{\rm  We can define extra three pseudo-spherical Lorentzian Darboux images along $\bm{\gamma}$:
\begin{align*}
(F)& \ \overline{\bm{D}}^T_n:I\lon H^2(-1);\ \overline{\bm{D}}^T_n(s) = \frac{\kappa_g\bn_{\bm{\gamma}}(s)+\kn\bb(s)}{\sqrt{\kg^2 - \kn^2}}\ \mbox{if}\ \kappa_g(s)^2>\kn^2, \\
(G)& \ \overline{\bm{D}}^S_n:I\lon S^2_1;\ \overline{\bm{D}}^S_n(s) = \frac{\kappa_g\bn_{\bm{\gamma}}(s)+\kn\bb(s)}{\sqrt{\kn^2 - \kg^2}}\ \mbox{if}\ \kappa_n(s)^2>\kg^2, \\
(H)& \ \overline{\bm{D}}^L_n:I\lon LC^*;\ \overline{\bm{D}}^L_n(s) = \frac{\kappa_g\bn_{\bm{\gamma}}(s)+\kn\bb(s)}{\sqrt{\kg^2 - \kn^2}}+\bm{t}(s)\ \mbox{if}\ \kappa_g(s)^2>\kn^2, \\
\end{align*}
Singularities and geometric meanings of (F) and (G) were investigated \cite{Sato}.
Moreover, we can easily obtain the similar results for (H), so that (H) was also essentially investigated in \cite{Sato}.
Therefore we omit the investigations on those three cases here.
}
\end{Rem}

\section{Singularities of pseudo-spherical Lorentzian Darboux images}
In this section we present a classification result of  
the singularities of pseudo-spherical Lorentzian Darboux images.
\par
We now introduce five invariants 
of $(M,\bgamma)$ as follows: 
\begin{align*}
(A)& \ \dtr (s) = \kn - \frac{\kg\tg'-\kg'\tg}{\kg^2 - \tg^2} \ \mbox{if}\ \kappa_g(s)^2>\tau_g(s)^2,\\           
(B)& \ \dsr (s)= \kn + \frac{\kg\tg'-\kg'\tg}{\tg^2 - \kg^2}\ \mbox{if}\ \tau_g(s)^2>\kappa_g(s)^2,\\           
(C)& \ \dlr (s)= \kn - \frac{\kg\tg'-\kg'\tg}{\kg^2 - \tg^2} + \sqrt{\kg^2 - \tg^2}\ \kappa_g(s)^2>\tau_g(s)^2,\\ 
(D)& \ \dso (s) = \kg + \frac{\kn\tg'-\kn'\tg}{\kn^2 + \tg^2}\ \mbox{if}\ (\kappa _n(s),\tau _g(s))\not= (0,0),\\          
(E)& \ \dlo (s)= \kg + \frac{\kn\tg'-\kn'\tg}{\kn^2 + \tg^2} + \sqrt{\kn^2 + \tg^2}\ \mbox{if}\ (\kappa _n(s),\tau _g(s))\not= (0,0). 
\end{align*}
We can classify the singular points of pseudo-spherical Lorentzian Darboux images
by using the above invariants.
\begin{Th}
Let $\gamma :I\lon M$ be a unit speed curve on a spacelike surface $M\subset \R^3_1$ such that $\Vert\boldsymbol{t}'(s)\Vert \neq 0$ and $\Vert \boldsymbol{b}'(s)\Vert \neq 0$.  
\par\noindent
{\rm (A)} Suppose that $\kappa_g ^2(s_0) > \tau_g ^2(s_0).$ Then we have the following assertions{\rm :}
		\par
		{\rm (1)} $\overline{\bm{D}}^T_r$ is non-singular at $s_0$ if and only if $\dtr(s_0) \neq 0.$
		\par
		{\rm (2)} The image of $\overline{\bm{D}}^T_r$ is locally diffeomorphic to the ordinary cusp $C$ at $s_0$ if and only if $\dtr(s_0) = 0$ and $(\dtr)'(s_0) \neq 0$	
		\par\noindent
{\rm (B)} Suppose that $\tau_g^2(s_0) > \kappa_g^2(s_0).$ Then we have the following assertions{\rm :}
		\par
		{\rm (1)} $\overline{\bm{D}}^S_r$ is non-singular at $s_0$ if and only if $\dsr(s_0) \neq 0.$
		\par
		{\rm (2)} The image of $\overline{\bm{D}}^S_r$ is locally diffeomorphic to the ordinary cusp $C$ at $s_0$ if and only if $\dsr(s_0) = 0$	and $(\dsr)'(s_0)\not= 0.$	

\par\noindent
{\rm (C)} Suppose that $\kappa_g^2(s_0) > \tau_g^2(s_0).$ Then we have the following assertions{\rm :}
		\par
		{\rm (1)} $\overline{\bm{D}}^L_r$ is non-singular at $s_0$ if and only if $\dlr(s_0)\not= 0.$
		\par
		{\rm (2)} The image of $\overline{\bm{D}}^L_r$ is locally diffeomorphic to the ordinary cusp $C$ at $s_0$ if and only if $\dlr(s_0) = 0$ and $(\dlr)'(s_0) \neq 0.$
		\par\noindent
{\rm (D)} Suppose that $(\kappa _n(s),\tau _g(s))\not= (0,0).$ Then we have the following assertions{\rm :}
		\par
		{\rm (1)} $\overline{\bm{D}}^S_o$ is non-singular at $s_0$ if and only if $\dso(s_0)\not= 0.$
		\par
		{\rm (2)} The image of $\overline{\bm{D}}^S_o$ is locally diffeomorphic to the ordinary cusp $C$ at $s_0$ if and only if $\dso(s_0) = 0$ and $(\dso)'(s_0) \neq 0.$
\par\noindent
{\rm (E)} Suppose that $(\kappa _n(s),\tau _g(s))\not= (0,0).$ Then we have the following assertions{\rm :}
		\par
		{\rm (1)} $\overline{\bm{D}}^L_o$ is non-singular at $s_0$ if and only if $\dlo(s_0)\not= 0.$
		\par
		{\rm (2)} The image of $\overline{\bm{D}}^L_o$ is locally diffeomorphic to the ordinary cusp $C$ at $s_0$ if and only if $\dlo(s_0) = 0$ and $(\dlo)'(s_0) \neq 0.$
		\end{Th}
	\par\noindent
		Here, $C=\{(x_1,x_2)\ |\ x_1^2=x_2^3\}$ is the {\it ordinary cusp} (or, the {\it semi-cubic parabola}).

\section{Legendrian dualities}
\par
We now review some properties of contact manifolds and Legendrian submanifolds.
Let $N$ be a $(2n+1)$-dimensional smooth manifold and $K$ be a 
tangent hyperplane field on $N$. 
Locally such a field is defined as the field of zeros of a $1$-form $\alpha .$
The tangent hyperplane field $K$ is {\it non-degenerate} if 
$\alpha \wedge (d\alpha )^n\not=0$ at any point of $N.$
We say that $(N,K)$ is a {\it contact manifold} if
$K$ is a non-degenerate hyperplane field.
In this case $K$ is called a {\it contact structure} and $\alpha $ is
a {\it contact form}.
Let $\phi :N\longrightarrow N'$ be a diffeomorphism between
contact manifolds $(N,K)$ and $(N',K').$
We say that $\phi $ is a {\it contact diffeomorphism} if $d\phi (K)=K'.$
Two contact manifolds $(N,K)$ and $(N',K')$ are
{\it contact diffeomorphic} if there exists a contact diffeomorphism
$\phi :N\longrightarrow N'.$
A submanifold $i:L\subset  N$ of  a contact manifold $(N,K)$ is said to be 
{\it Legendrian} if ${\rm dim}\ L=n$ and $di_x(T_xL)\subset K_{i(x)}$ at any
$x\in L.$  We say that a smooth fiber bundle $\pi :E\lon M$ is
called a {\it Legendrian fibration} if its total space $E$ is
furnished with a contact structure and its fibers are Legendrian submanifolds.
Let $\pi :E\lon M$ be a Legendrian fibration.
For a Legendrian submanifold $i:L\subset E,$ $\pi\circ i:L\lon M$
is called a {\it Legendrian map.}  The image of the Legendrian map
$\pi\circ i$ is called a {\it wavefront set} of $i$ which is
denoted by $W(L).$
For any $z\in E,$ it is known that there is a local coordinate
system $(x,p,y)=(x_1,\dotsc ,x_m,p_1,\dotsc ,p_m,y)$ around $z$ such that
$
\pi (x,p,y)=(x,y)
$
and the contact structure is given by the 1-form
$\alpha =dy-\sum_{i=1}^m p_idx_i
$
(cf. \cite{Arnold1}, 20.3).
\par
In \cite{Izu2} we have shown the basic duality theorem
which is the fundamental tool for the study of
spacelike hypersurfaces in Lorentz-Minkowski pseudo-spheres.
We consider the following four double fibrations:
\par\noindent
(1) (a) $H^2(-1)\times S^2_1\supset \Delta _1
=\{(\bv ,\bw )\ |\ \langle \bv ,\bw \rangle =0\ \}$,
\par
(b) $\pi _{11}:\Delta _1\lon H^2(-1)$,$\pi _{12}:\Delta _1\lon S^2_1$,
\par
(c) $\theta _{11}=\langle d\bv ,\bw\rangle |\Delta _1$,
$\theta _{12}=\langle \bv ,d\bw\rangle |\Delta _1$.
\bigskip
\par\noindent
(2) (a) $H^2(-1)\times LC^*\supset \Delta _2=
\{(\bv ,\bw )\ |\ \langle \bv ,\bw \rangle =-1\ \}$,
\par
(b) $\pi _{21}:\Delta _2\lon H^2(-1)$,$\pi _{22}:\Delta _2\lon LC^*$,
\par
(c) $\theta _{21}=\langle d\bv ,\bw\rangle |\Delta _2$,
$\theta _{22}=\langle \bv ,d\bw\rangle |\Delta _2$.
\bigskip
\par\noindent
(3) (a) $LC^*\times S^2_1\supset \Delta _3=
\{(\bv ,\bw )\ |\ \langle \bv ,\bw \rangle =1\ \}$,
\par
(b) $\pi _{31}:\Delta _3\lon LC^*$,$\pi _{32}:\Delta _3\lon S^2_1$,
\par
(c) $\theta _{31}=\langle d\bv ,\bw\rangle |\Delta _3$,
$\theta _{32}=\langle \bv ,d\bw\rangle |\Delta _3$.
\bigskip
\par\noindent
(4) (a) $LC^*\times LC^*\supset \Delta _4=
\{(\bv ,\bw )\ |\ \langle \bv ,\bw \rangle =-2\ \}$,
\par
(b) $\pi _{41}:\Delta _4\lon LC^*$,$\pi _{42}:\Delta _4\lon LC^*$,
\par
(c) $\theta _{41}=\langle d\bv ,\bw\rangle |\Delta _4$,
$\theta _{42}=\langle \bv ,d\bw\rangle |\Delta _4$.
\par
Here, $\pi _{i1}(\bv,\bw)=\bv$, $\pi _{i2}(\bv ,\bw)=\bw$,
$\langle d\bv,\bw\rangle =-w_0dv_0+\sum_{i=1}^2 w_idv_i$
and
$\langle \bv,d\bw \rangle =-v_0dw_0+\sum_{i=1}^2 v_idw_i.$
\par
We remark that 
$\theta _{i1}^{-1}(0)$ and $\theta _{i2}^{-1}(0)$
define the same tangent hyperplane field over $\Delta _i$ which is denoted by $K_i.$
The basic duality theorem is the following theorem \cite{Izu2}:
\begin{Th}
 \label{thm:basic-dual}
 With the same notations as the previous paragraph,
 each $(\Delta _i,K_i)$ $(i=1,2,3,4)$ is a contact manifold and
 both of $\pi _{ij}$ $(j=1,2)$ are Legendrian fibrations.
 Moreover those contact manifolds are contact diffeomorphic each other.
\end{Th}
Moreover, we have the following extra double fibration:
\begin{enumerate}
\item[(5)]
   \begin{enumerate}
    \item[(a)]$S^2_1 \times S^2_1 \supset \Delta_5 = \lbrace(\bv,\bw) \ | \ \langle \bv,\bw \rangle = 0\rbrace,$
    \item[(b)]$\pi_{51} : \Delta_5 \longrightarrow S^2_1,\pi_{52} : \Delta_1 \longrightarrow S^2_1,$
    \item[(c)]$\theta_{51} = \langle d\bv,\bw \rangle | \Delta_5,\theta_{52} = \langle \bv,d\bw \rangle | \Delta_5.$         
   \end{enumerate}
\end{enumerate}
It is shown in \cite{CIZ09} that $(\Delta _5, K_5)$ is a contact manifold such that $\pi _{5j}:\Delta_5 \longrightarrow S^2_1$, $j=1,2$, are
Legendrian fibrations which is not contact diffeomorphic to any other $(\Delta _i,K_i),$ $i=1,2,3,4.$
Given a Legendrian submanifold $i:L\to \Delta_i$, $i=1,2,3,4,5$,
We say that $\pi_{i1}(i(L))$ is the {\it $\Delta _i$-dual
of $\pi_{i2}(i(L))$}  and vice-versa.

Then we have the following duality theorem.
\begin{Th} Let $\gamma :I\lon M$ be a unit speed curve on a spacelike surface $M\subset \R^3_1$ such that $\Vert\boldsymbol{t}'(s)\Vert \neq 0$ and $\Vert \boldsymbol{b}'(s)\Vert \neq 0$.  
\begin{enumerate}
\item[{\rm (1)}] If $(\kn,\tg) \neq (0,0)$, then $\bn_{\bm{\gamma}}$ is a $\Delta_1$-dual of $\overline{\bm{D}}^S_o$. 
 \item[{\rm (2)}] If $(\kn,\tg) \neq (0,0)$, then $\bn_{\bm{\gamma}}$ is a $\Delta_2$-dual of $\overline{\bm{D}}^L_o$.
 \item[{\rm (3)}] If $\kg^2 > \tg^2$, then $\bb$ is a $\Delta_1$-dual of $\overline{\bm{D}}^T_r$. 
  \item[{\rm (4)}] If $\kg^2 > \tg^2$, then $\bb$ is a $\Delta_3$-dual of $\overline{\bm{D}}^L_r$.
 \item[{\rm (5)}] If $\tg^2 > \kg^2$, then $\bb$ is a $\Delta_5$-dual of $\overline{\bm{D}}^S_r$. 
\end{enumerate}
\end{Th}
\demo We can show that (1) holds as follows:
\begin{enumerate}
\item[{\rm (1)}]
We define a mapping 
${\cal L}_1 : I  \longrightarrow \Delta_1$ by ${\cal L}_1(s) = (\bn_{\bm{\gamma}}(s),\overline{\bm{D}}^S_o(s))$.
Then  we have $\langle\bn_{\bm{\gamma}}(s),\overline{\bm{D}}^S_o(s)\rangle = 0$ and  ${\cal L}^*_1\theta_{11} = \langle\boldsymbol{n}'_{\gamma}(s),\overline{\bm{D}}^S_o(s)\rangle = 0$. Thus ${\cal L}_1$ is an isotropic mapping, so that
$\bn_{\bm{\gamma}}$ is a $\Delta_1$-dual of $\overline{\bm{D}}^S_o$.
\end{enumerate}
Then we define mappings
\begin{eqnarray*}
{\cal L}_2 : I  \longrightarrow \Delta_2 &;& {\cal L}_2(s) = (\bn_{\bm{\gamma}}(s),\overline{\bm{D}}^L_o(s)), \\
{\cal L}_3 : I  \longrightarrow \Delta_1 &;& {\cal L}_3(s) = (\bb(s),\overline{\bm{D}}^T_r(s)), \\
{\cal L}_4 : I \longrightarrow \Delta_3 &;& {\cal L}_3(s) = (\bb(s),\overline{\bm{D}}^L_r(s))), \\
{\cal L}_5 : I \longrightarrow \Delta_5 &;& {\cal L}_5(s) = (\bb(s),\overline{\bm{D}}^S_r(s)).
\end{eqnarray*}
Then we can show that $\mathcal{L}_i$ $(i=2,3,4,5)$ are isotropic mappings.
This means that (2), (3), (4) and  (5) hold.
\enD

\section{Height functions}
 We now introduce five families of functions on $\bm{\gamma}:I \longrightarrow M \subset \R^3_1$
 with $\Vert\boldsymbol{t}'(s)\Vert \neq 0 ,\Vert \boldsymbol{b}'(s)\Vert \neq 0$ as follows:
\begin{align*}
H^T_r :& \ I \times H^2_{+}(-1) \longrightarrow \R \ ; \ (s,\bv) \longmapsto \langle \bb(s),\bv \rangle, \\
H^S_r :& \ I \times S^2_1 \longrightarrow \R \ ; \ (s,\bv) \longmapsto \langle \bb(s),\bv \rangle,\\  
H^L_r :& \ I \times LC^* \longrightarrow \R \ ; \ (s,\bv) \longmapsto \langle \bb(s),\bv \rangle - 1,\\ 
H^S_o :& \ I \times S^2_1 \longrightarrow \R \ ; \ (s,\bv) \longmapsto \langle \bn_{\bm{\gamma}}(s),\bv \rangle,\\
H^L_o :& \ I \times LC^* \longrightarrow \R \ ; \ (s,\bv) \longmapsto \langle \bn_{\bm{\gamma}}(s),\bv \rangle + 1. 
\end{align*}
For any $\bm{v},$ we define $h^T_{r,\bm{v}}(s)= H^T_r(s,\bm{v}),h^S_{r,\bm{v}}(s)= H^S_r(s,\bm{v}),h^L_{r,\bm{v}}(s)= H^L_r(s,\bm{v}),h^S_{o,\bm{v}}(s)= H^S_o(s,\bm{v}),h^L_{o,\bm{v}}(s)= H^L_o(s,\bm{v})$
Then we have the following proposition.
\begin{Pro} Let $\bm{\gamma}:I \longrightarrow M $ be a unit speed curve on a spacelike surface $M\subset \R^3_1$ such that
$\Vert\boldsymbol{t}'(s)\Vert \neq 0 ,\Vert \boldsymbol{b}'(s)\Vert \neq 0$. Then we have the following\/{\rm :}
\begin{enumerate}
  \item[{\rm (A)}] For any $(s,\bm{v}) \in I \times H^2(-1)$, we have the following\/{\rm :}
 \begin{enumerate}
       \item[{\rm (1)}] $h^T_{r,\bm{v}}(s) = 0$ if and only if there exist $\lambda,\mu \in \R$ with $-\lambda^2 + \mu^2 = 1$ such that     
       \[ \bm{v} = \lambda\bt(s) + \mu\bn_{\bm{\gamma}}(s), \]
       \item[{\rm (2)}] $h^T_{r,\bm{v}}(s) =(h^T_{r,\bm{v}})'(s) = 0$ if and only if $\kg^2 > \tg^2$ and
                  $\bm{v} = \pm\overline{\bm{D}}^T_r(s),$
       \item[{\rm (3)}] $h^T_{r,\bm{v}}(s) =(h^T_{r,\bm{v}})'(s) = (h^T_{r,\bm{v}})''(s) = 0$ if and only if $ \kg^2 > \tg^2,\ \dtr (s) = 0$ and
                 $ \bv = \pm\overline{\bm{D}}^T_r(s),$
       \item[{\rm (4)}] $h^T_{r,\bm{v}}(s) =(h^T_{r,\bm{v}})'(s) = (h^T_{r,\bm{v}})''(s) = (h^T_{r,\bm{v}})'''(s) = 0$ if and only if $ \kg^2 > \tg^2,$ $\dtr(s) = 0$, $(\dtr)'(s) = 0$ and
       $\bv = \pm\overline{\bm{D}}^T_r(s).$ 

     \end{enumerate}
\item[{\rm (B)}] For any $(s,\bm{v}) \in I \times S^2_1$, we have the following\/{\rm :}
     \begin{enumerate}
      \item[{\rm (1)}] $h^S_{r,\bm{v}}(s) = 0$ if and only if there exist $\lambda,\mu \in \R$ with $-\lambda^2 + \mu^2 = -1$ such that    
       \[ \bm{v} = \lambda\bt(s) + \mu\bn_{\bm{\gamma}}(s), \]
       \item[{\rm (2)}] $h^S_{r,\bm{v}}(s) =(h^S_{r,\bm{v}})'(s) = 0$ if and only if $\tg^2 > \kg^2$ and 
                  $\bm{v} = \pm\overline{\bm{D}}^S_r(s),$ 
       \item[{\rm (3)}] $h^S_{r,\bm{v}}(s) =(h^S_{r,\bm{v}})'(s) = (h^S_{r,\bm{v}})''(s) = 0$ if and only if $ \tg^2 > \kg^2,$ $\dsr(s) = 0$ and
       $\bm{v} = \pm\overline{\bm{D}}^S_r(s),$ 
       \item[{\rm (4)}] $h^S_{r,\bm{v}}(s) =(h^S_{r,\bm{v}})'(s) = (h^S_{r,\bm{v}})''(s) = (h^S_{r,\bm{v}})'''(s) = 0$ if and only if $ \tg^2 > \kg^2,$ $\dsr(s) = 0,$ $ (\dsr)'(s) = 0$
       and $\bm{v} = \pm\overline{\bm{D}}^S_r(s).$ 
     \end{enumerate} 
      \item[{\rm (C)}] For any $(s,\bm{v}) \in I \times LC^*$, we have the following\/{\rm :}
     \begin{enumerate}
     \item[{\rm (1)}] $h^L_{r,\bm{v}}(s) = 0$ if and only if there exist $\lambda,\mu \in \R$ with $\lambda^2 - \mu^2 = -1$ such that     
       \[ \bm{v} = \lambda\bt(s) + \mu\bn_{\bm{\gamma}}(s) +\bb(s),\]
       \item[{\rm (2)}] $h^L_{r,\bm{v}}(s) =(h^L_{r,\bm{v}})'(s) = 0$ if and only if $\kg^2 > \tg^2$ and
                  $\bm{v} = \pm\overline{\bm{D}}^L_r(s),$
       \item[{\rm (3)}] $h^L_{r,\bm{v}}(s) =(h^L_{r,\bm{v}})'(s) = (h^L_{r,\bm{v}})''(s) = 0$ if and only if $ \kg^2 > \tg^2,$ $\dlr(s) = 0$ and
                 $\bm{v} = \pm\overline{\bm{D}}^L_r(s),$
       \item[{\rm (4)}] $h^L_{r,\bm{v}}(s) =(h^L_{r,\bm{v}})'(s) = (h^L_{r,\bm{v}})''(s) = (h^L_{r,\bm{v}})'''(s) = 0$ if and only if  $\kg^2 > \tg^2,$ $\dlr (s) = 0,$ $ (\dlr)'(s) = 0$ and 
       $\bm{v} = \pm\overline{\bm{D}}^L_r(s).$
     \end{enumerate}
    \item[{\rm (D)}] Suppose that $(\kappa _n(s),\tau _g(s))\not= (0,0).$ For any $(s,\bm{v}) \in I \times S^2_1$, we have the following\/{\rm :}
     \begin{enumerate}
     \item[{\rm (1)}] $h^S_{o,\bm{v}}(s) = 0$ if and only if there exist $\lambda,\mu \in \R$ with $\lambda^2 + \mu^2 = 1$ such that   
       \[ \bv = \lambda\bt(s) + \mu\bb(s), \]
       \item[{\rm (2)}] $h^S_{o,\bm{v}}(s) =(h^S_{o,\bm{v}})'(s) = 0$ if and only if 
                   $\bm{v} = \pm\overline{\bm{D}}^S_o(s),$ 
       \item[{\rm (3)}] $h^S_{o,\bm{v}}(s) =(h^S_{o,\bm{v}})'(s) = (h^S_{o,\bm{v}})''(s) = 0$ if and only if  $\dso(s) = 0$ and
                 $\bm{v} = \pm\overline{\bm{D}}^S_o(s),$
       \item[{\rm (4)}] $h^S_{o,\bm{v}}(s) =(h^S_{o,\bm{v}})'(s) = (h^S_{o,\bm{v}})''(s) = (h^S_{o,\bm{v}})'''(s) = 0$ if and only if  $\dso(s) = 0,$ $ (\dso)'(s) = 0$ and
       $\bm{v} = \pm\overline{\bm{D}}^S_o(s).$
         \end{enumerate}
   \item[{\rm (E)}] Suppose that $(\kappa _n(s),\tau _g(s))\not= (0,0).$ For any $(s,\bm{v}) \in I \times LC^*$, we have the following\/{\rm :}
     \begin{enumerate}
     \item[{\rm (1)}] $h^L_{o,\bm{v}}(s) = 0$ if and only if there exist $\lambda,\mu \in \R$ with $\lambda^2 + \mu^2 = 1$ such that     
       \[ \bm{v} = \lambda\bt(s) + \mu\bb(s) + \bn_{\bm{\gamma}}(s), \]
       \item[{\rm (2)}] $h^L_{o,\bm{v}}(s) =(h^L_{o,\bm{v}})'(s) = 0$ if and only if 
        $\bm{v} = \pm\overline{\bm{D}}^L_o(s),$
       \item[{\rm (3)}] $h^L_{o,\bm{v}}(s) =(h^L_{o,\bm{v}})'(s) = (h^L_{o,\bm{v}})''(s) = 0$ if and only if  $\dlo(s) = 0$ and
                 $\bm{v} = \pm\overline{\bm{D}}^L_o(s),$
       \item[{\rm (4)}] $h^L_{o,\bm{v}}(s) =(h^L_{o,\bm{v}})'(s) = (h^L_{o,\bm{v}})''(s) = (h^L_{o,\bm{v}})'''(s) = 0$ if and only if $\dlo(s) = 0,$ $ (\dlo)'(s) = 0$ and
       $\bm{v} = \pm\overline{\bm{D}}^L_o(s).$
     \end{enumerate}
     \end{enumerate}
\end{Pro}
\demo
We remark that $\Vert\boldsymbol{t}'(s)\Vert \neq 0 ,\ \Vert \boldsymbol{b}'(s)\Vert \neq 0$ if and only if
$-\kn^2+\kg^2\not= 0,\ \kg^2-\tg^2\not=0. $
\begin{enumerate}
\item[(A)]By straight forward calculations, we have the following :
  \begin{align*}
   h^T_{r,\bm{v}} =& \langle\bm{b},\bm{v} \rangle,\\ 
   (h^T_{r,\bm{v}})' =& \langle -\kgs \bm{t}+\tgs\bm{n}_{\bm{\gamma}},\bm{v} \rangle,\\
   (h^T_{r,\bm{v}})'' =& \langle(-\kgs'+\kns\tgs)\bm{t}-(\kgs^2-\tgs^2)\bm{b}+(\tgs'-\kgs\kns)\bm{n}_{\bm{\gamma}},\bm{v} \rangle,\\
   (h^T_{r,\bm{v}})''' =& \langle(-\kgs''+\kns'\tgs+2\kns\tgs'+\kgs^3-\kgs\tgs^2-\kns^2\kgs)\bm{t}\\
                     &+(-3\kgs\kgs'+3\tgs\tgs')\bm{b}+(\tgs''-\kgs\kns'-2\kns\kgs'+\tgs^3+\kns^2\tgs-\kgs^2\tgs)\bm{n}_{\bm{\gamma}},\bm{v}\rangle.
  \end{align*}
   \end{enumerate}
Since $\{\bm{t}(s), \bm{n}_{\bm{\gamma}}(s),\bm{b}(s)\}$ is a pseudo-orthonormal frame of $\R^3_1$ along $\bm{\gamma},$ we have
$\bm{v}=\lambda \bm{t}(s)+\mu \bm{n}_{\bm{\gamma}}(s)+ \eta \bm{b}(s).$
\begin{enumerate}
   \item[1.] Since $h^T_{r,\bm{v}} = 0$, $\eta =0,$ so that we have $\bm{v}=\lambda \bm{t}(s)+\mu \bm{n}_{\bm{\gamma}}(s).$
   Here, $\bm{v}$ is timelike. Then we have $\mu^2>\lambda^2$. Thus, $\mu\not= 0.$ This completes the proof of assertion (A),(1).
   \item[2.] Moreover, $(h^T_{r,\bv})' = 0$ implies $-\lambda\kgs-\mu\tgs = 0$. Therefore, we have
   $\kgs \bm{v}=\kgs\lambda \bm{t}+\kgs\mu\bm{n}_{\bm{\gamma}}=-\mu(\tgs\bm{t}-\kgs\bm{n}_{\bm{\gamma}}).$
   Thus we have $-\kgs^2=\mu^2(\tgs^2-\kgs^2),$ so that $\tgs^2\leq \kgs^2.$
   Since $\kg^2-\tg^2\not=0,$ we have  $\tgs^2< \kgs^2.$
   It follows that
   \[\bm{v} =\pm\frac{\tgs\ts-\kgs\ns}{\sqrt{\kgs^2 - \tgs^2}}=\pm\overline{\bm{D}}^T_r.\]
   \item[3.] If we add extra condition $(h^T_{r,\bm{v}})'' = 0$, then we have
   \[ \kns - \frac{\kgs\tgs'-\kgs'\tgs}{\kgs^2 - \tgs^2} = 0. \]
   Thus we have $\dtr=0.$
   \item[4.]Moreover, if we consider one more condition $(h^T_{r,\bm{v}})'''=0$, then we have
   \[ \kns'(\kgs^2 - \tgs^2)+2(\kgs\kgs'-\tgs\tgs')(\dtr + \frac{\kgs\tgs'-\kgs'\tgs}{\kgs^2 - \tgs^2}) - (\kgs\tgs'' - \kgs''\tgs) = 0 .\]
Since we have
   \[(\dtr)' =  \kns' +\frac{2(\kgs\kgs'-\tgs\tgs')(\kgs\tgs'-\kgs'\tgs)}{(\kgs^2 - \tgs^2)^2} -\frac{(\kgs\tgs'' - \kgs''\tgs)}{(\kgs^2 - \tgs^2)},\]
   $(\dtr)'=0$ with the condition $\dtr=0$.
   \end{enumerate}
\par
For other cases (B), (C), (D) and (E), we
have the similar calculations to case (A) for the derivatives of $h^S_{r,\bm{v}},h^L_{r,\bm{v}},h^S_{o,\bm{v}}$ and $h^L_{o,\bm{v}}$, respectively.
We omit the details here.
\enD
\section{Proof of Theorem 3.1}
In this section we give a proof of Theorem 3.1.
In order to prove Theorem 3.1, we use some general results on the singularity theory for families of function germs.
Detailed descriptions are found in the book\cite{Bru-Gib}.
Let $F:({\R}\times{\R}^r,(s_0,x_0))\longrightarrow {\R}$ be a function germ.
We call $F$ an {\it $r$-parameter unfolding\/} of $f$, where $f(s)=F_{x_0}(s,x_0).$
We say that $f$ has an {\it $A_k$-singularity\/} at $s_0$ if $f^{(p)}(s_0)=0$ for all $1\leq p\leq k$,
and $f^{(k+1)}(s_0)\ne 0.$
We also say that $f$ has an {\it $A_{\ge k}$-singularity\/} at $s_0$ if
$f^{(p)}(s_0)=0$ for all $1\leq p\leq k.$ 
Let $F$ be an unfolding of $f$ and $f(s)$ has an $A_k$-singularity $(k\geq 1)$ at $s_0.$
We denote the $(k-1)$-jet of the partial derivative
$\frac{\partial F}{\partial x_i}$ at $s_0$ by 
$j^{(k-1)}(\frac{\partial F}{\partial x_i}(s,x_0))(s_0)=\sum_{j=0}^{k-1} \alpha_{ji}(s-s_0)^j$ 
for $i=1,\dots ,r$. Then
$F$ is called an {\it $\mathcal{R}$-versal unfolding \/} if the $k\times{r}$ 
matrix of coefficients $(\alpha _{ji})_{j=0,\dots ,k-1;i=1,\dots ,r}$
has rank $k$ $(k\leq {r}).$
We introduce an important set concerning the unfoldings relative to the above notions. 
The {\it discriminant set} of $F$ is the set
$${\mathcal D}_F=\{x\in {\R}^r|{\rm  there\ exists }\ s \ {\rm with }\ F=
\frac{\partial F}{\partial s}=0 \ {\rm at }\ (s,x)\}.$$
Then we have the following classification (cf., \cite{Bru-Gib}).
\begin{Th} Let $F:({\R}\times{\R}^r,(s_0,x_0))\longrightarrow {\R}$
be an $r$-parameter unfolding of $f(s)$ which has the $A_2$ singularity at $s_0$.
If $F$ is an $\mathcal{R}$-versal unfolding, then ${\mathcal D}_F$ is locally diffeomorphic to  
$C\times {\R}^{r-1}$.
\end{Th}
Here, $C=\{(x_1,x_2)\ |\ x_1=t^2, x_2=t^3\}$ is the {\it ordinary cusp} (i.e. the {\it semi-cubic parabola}).

We now consider that 
$H^T_r,
H^S_r, 
H^L_r,H^T_o,$ and
$H^L_o$ are unfoldings of
$h^T_{r,\bm{v}}(s),h^S_{r,\bm{v}}(s),h^L_{r,\bm{v}}(s),h^T_{o,\bm{v}}(s),$ and $h^L_{o,\bm{v}}(s)$ for any $\bm{v},$ respectively.

\begin{Pro}  Let $\bm{\gamma}:I \longrightarrow M $ be a unit speed curve on a spacelike surface $M\subset \R^3_1$ such that
$\Vert\boldsymbol{t}'(s)\Vert \neq 0 ,\Vert \boldsymbol{b}'(s)\Vert \neq 0$. Then we have the following\/{\rm :}
\begin{enumerate}

  \item[{\rm (A)}] if $h^T_{r,\bv}$ has the $A_2$-singularity at $s_0$, then $H^T_r$ is an $\mathcal{R}$-versal unfolding of $h^T_{r,\bv}$,
  \item[{\rm (B)}] if $h^S_{r,\bv}$ has the $A_2$-singularity at $s_0$, then $H^S_r$ is an $\mathcal{R}$-versal unfolding of $h^S_{r,\bv}$,
  \item[{\rm (C)}] if $h^L_{r,\bv}$ has the $A_2$-singularity at $s_0$, then $H^L_r$ is an $\mathcal{R}$-versal unfolding of $h^L_{r,\bv}$,
  \item[{\rm (D)}] if $h^S_{o,\bv}$ has the $A_2$-singularity at $s_0$, then $H^S_o$ is an $\mathcal{R}$-versal unfolding of $h^S_{o,\bv}$,
  \item[{\rm (E)}] if $h^L_{o,\bv}$ has the $A_2$-singularity at $s_0$, then $H^L_o$ is an $\mathcal{R}$-versal unfolding of $h^L_{o,\bv}$. 
\end{enumerate}
\end{Pro}
\demo Here, we only give the proof for (A). Other cases are similar to case (A).
\par\noindent
(A) We denote that $\bv = (\sqrt{x_1^2 + x_2^2 + 1},x_1,x_2) \in H^2_+(-1),\bb = (b_0(s),b_1(s),b_2(s))$. Then we have
\[H^T_r(s,\bv) = - b_0(s)\sqrt{x_1^2 + x_2^2 + 1}+ b_1(s)x_1 + b_2(s)x_2 \]
and
 \[
\frac{\partial{H^T_r}}{\partial{x_1}}(s,\bv) = -b_0(s)\frac{x_1}{\sqrt{x_1^2 + x_2^2 + 1}} + b_1(s) \ , \ \frac{\partial{H^T_r}}{\partial{x_1}}(s,\bv) = -b_0(s)\frac{x_2}{\sqrt{x_1^2 + x_2^2 + 1}} +b_2(s).
\]
Therefore the $2$-jet of $H^T_r(s,\bv)$ are
\begin{eqnarray*}
j^2\frac{\partial{H^T_r}}{\partial{x_1}}(s_0,\bv) &=& \left( -b_0(s_0)\frac{x_1}{\sqrt{x_1^2 + x_2^2 + 1}} + b_1(s_0) \right) \\
&{}&+ \left(  -b'_0(s_0)\frac{x_1}{\sqrt{x_1^2 + x_2^2 + 1}}+b'_1(s_0)  \right)(s- s_0),\\
j^2\frac{\partial{H^T_r}}{\partial{x_2}}(s_0,\bv) &=& \left( -b_0(s_0)\frac{x_2}{\sqrt{x_1^2 + x_2^2 + 1}} 
+ b_2(s_0) \right) \\
&{}& + \left(  -b'_0(s_0)\frac{x_2}{\sqrt{x_1^2 + x_2^2 + 1}}+b'_2(s_0)  \right)(s- s_0).
\end{eqnarray*}
We consider the following matrix:
\[A=
\begin{pmatrix}
-b_0(s_0)\frac{x_1}{\sqrt{x_1^2 + x_2^2 + 1}} + b_1(s_0) & -b_0(s_0)\frac{x_2}{\sqrt{x_1^2 + x_2^2 + 1}} 
+ b_2(s_0)\\
 -b'_0(s_0)\frac{x_1}{\sqrt{x_1^2 + x_2^2 + 1}}+b'_1(s_0)  & -b'_0(s_0)\frac{x_2}{\sqrt{x_1^2 + x_2^2 + 1}}+b'_2(s_0)
\end{pmatrix}
\]
The determinant of $A$ is
\begin{eqnarray*}
\mbox{det}A &=&\frac{1}{\sqrt{x_1^2 +x_2^2 +1}}\left( x_1(b'_0b_2 - b_0b'_2) + x_2(b_0b'_1 - b'_0b_1) + \sqrt{x_1^2 + x_2^2 + 1}(b'_2b_1-b_2b'_1) \right)\\
            &=&\frac{1}{\sqrt{x_1^2 +x_2^2 +1}}\left\langle(b_2b'_1 - b_1b'_2,b_2b'_0 - b'_2b_0,b_0b'_1-b'_0b_1),(\sqrt{x_1^2 + x_2^2 + 1},x_1,x_2)\right\rangle\\   
            &=&\frac{1}{\sqrt{x_1^2 + x_2^2 + 1}}\langle(\boldsymbol{b} \wedge \boldsymbol{b}'),\boldsymbol{v}\rangle\\
            &=&\frac{1}{\sqrt{x_1^2 + x_2^2 + 1}}\vert\boldsymbol{b} \  \ \boldsymbol{b}' \ \  \boldsymbol{v}\vert      
\end{eqnarray*}
By Proposition 5.1, if $h^T_{r,\bm{v}}$ has the $A_2$-singularity at $s_0$,  then $\bm{v}=\pm \overline{\bm{D}}^T_r(s)$
and $\bm{b}'(s)=-\kappa _g(s)\bm{t}(s)+\tau _g(s)\bm{n}_{\bm{\gamma}}(s)$, so that $\lbrace\boldsymbol{b},\boldsymbol{b}',\boldsymbol{v}\rbrace$ is linearly independent. Therefore, ${\rm rank}\, A=2.$  This means that  $H^T_r$ is an $\mathcal{R}$-versal unfolding of $h^T_{r,\bm{v}}.$
\enD

We define three vector fields respectively defined as normalizations of $\bm{t}', \bn_\gamma', \bb'$ as follows:
\[\boldsymbol{T}_{\boldsymbol{t}}(s) = \frac{\kn\bn_{\bm{\gamma}}(s)+\kg\bb(s)}{\sqrt{\kg^2 - \kn^2}},
\boldsymbol{T}_{\boldsymbol{n}_{\gamma}}(s) = \frac{\kn\bt(s)+\tg\bb(s)}{\sqrt{\kg^2 + \tg^2}},
\boldsymbol{T}_{\boldsymbol{b}}(s) = \frac{-\kg\bt(s)+\tg\bn_{\bm{\gamma}}(s)}{\sqrt{\kg^2 - \tg^2}}.
 \]

We can prove Theorem 3.1.
\smallskip 
\par\noindent
{\it Proof of Theorem 3.1.}
Here we only give the proof for (A) again.
\par\noindent
(A) (1) By a straight forward calculation $\overline{\bm{D}}^T_r(s)$, we have
 \[\left(\overline{\bm{D}}^T_r(s)\right)' = \dtr(s)\boldsymbol{T}_{\boldsymbol{b}}(s), \]
 so that $\overline{\bm{D}}^T_r(s)$ is non-singular at $s=s_0$ if and only if $\dtr(s_0)\not= 0.$
\par
(2) By Proposition 5.1, $h^T_{r,\bv}$ is an $A_2$-singularity if and only if
$\dtr(s_0) = 0$, $(\dtr)'(s_0) \neq 0$ and $\bm{v}=\overline{\bm{D}}^T_r(s).$
By Proposition 6.2, $H^T_r$ is an $\mathcal{R}$-versal unfolding of $h^T_{r,\bm{v}}$.
By Proposition 5.1, the image of $\overline{\bm{D}}^T_r$ is the discriminant set of $H^T_r.$
By Theorem 6.1, the discriminant of $H^T_r$ is locally diffeomorphic to the cusp $C.$
\enD
\section{Invariants of curves on surfaces}
In this section we consider geometric meanings of the invariants $\dtr,\dsr,\dlr,\dso,\dlo$.
In particular what can we say about the original curve on the surface when each invariant is vanishing.
For the purpose, we consider cylinders in Lorentz-Minkowski space.
A (generalized) {\it cylinder} in $\R^3_1$ is a ruled surface with a constant director.
It is parametrized by $F(t,u)=\bsigma (t)+u\bv,$ where $\bsigma$ is a smooth curve and $\bv$ is a 
non-zero vector. The vector $\bv$ is called the {\it director}.
We say that $F$ is a {\it spacelike cylinder}, a {\it timelike cylinder} and a {\it lightlike cylinder}
if the director $\bv$ is spacelike, timelike and lightlike, respectively.
Let $M$ be a surface and $N$ be a cylinder in $\R^3_1.$ We say that $N$ is a {\it pseudo-normal cylinder} of $M$ if
$M\cap N\not=\emptyset$ and 
$T_pN$ contains the pseudo-normal vector $\bn ( p)$ at any $p\in M\cap N.$
In this case $M$ and $N$ transversally intersect, so that $M\cap N$ is a regular curve $C$.
We call $C$ a {\it slice of $M$ with a pseudo-normal cylinder} of $M.$
Moreover, we call $N$ a {\it pseudo-normal spacelike cylinder} if the director of $N$ is spacelike and
a {\it pseudo-normal timelike cylinder} if the director of $N$ is timelike, respectively.
We remark that the director of $N$ is not lightlike.
If $N$ is locally parametrized by $F(t,u)=\bsigma (t)+u\bv,$ then we have
\[
\frac{\partial F}{\partial t}(t,u)=\bsigma '(t)\ \mbox{and}\ \frac{\partial F}{\partial u}(t,u)=\bv,
\]
so that the pseudo-normal to $N$ is given by
\[
\frac{\partial F}{\partial t}(t,u)\times \frac{\partial F}{\partial u}(t,u) =\bsigma '(t)\times \bv.
\]
If $C$ is parametrized by $\bgamma (s)$, where $s$ is the arc-length parameter of $\bgamma ,$
then $N$ can be parameterized by $F(s,u)=\bgamma (s)+u\bv$ at least locally.
Since $N$ has been given an orientation by $F,$ the unit normal vector of $N$ along $C$ is $\bb(s).$
In particular, $\langle \bb(s),\bv\rangle =0.$
\par
On the other hand, $N$ is called a {\it osculating cylinder} if the tangent planes of $M$ and $N$ coincide at
any point of $M\cap N.$
In this case $C=M\cap N$ is called a {\it slice of $M$ with an osculating cylinder} of $M.$
We remark that the director of the osculating cylinder is always spacelike.
If $N$ is locally parametrized by $F(t,u) = \boldsymbol{\sigma}(t) + u\bv,$ then 
the unit normal vector of $N$ along $C$ is $\bn_{\bm{\gamma}}$ and $\langle \bn ,\bv \rangle = 0$
for a parmetrization $\bm{\gamma}$ of $C.$
\par
We call $N$ a {\it hyperbolic lightlike cylinder} if $M \cap N \neq \emptyset$ and
$\langle \boldsymbol{n}(p),\bv \rangle = -1$ at any point $p\in M\cap N,$ where $\bv$ is the lightlike director of $N.$
In this case, $N$ is transversely intersect with $M$, so that $C=M\cap N$ is a regular curve.
We call $C$ a {\it slice of $M$ with a hyperbolic lightlike cylinder}.
We also call $N$ a {\it de Sitter lightlike cylinder} if $M \cap N \neq \emptyset$ and
$\langle \boldsymbol{b}(p),\bv \rangle = 1$ at any point $p\in M\cap N,$ where $\bv$ is the lightlike director of $N.$
In this case, $N$ is transversely intersect with $M$, so that $C=M\cap N$ is a regular curve.
We call $C$ a {\it slice of $M$ with a de Sitter lightlike cylinder}.
For the both cases in the above, $\bv$ is lightlike.

Then we have the following theorem.
\begin{Th}\label{theorem:invariant}
Let $\gamma:I \longrightarrow M$ be a unit speed curve on a spacelike surface $M\subset \R^3_1$ such that
$\Vert\boldsymbol{t}'(s)\Vert \neq 0 ,\Vert \boldsymbol{b}'(s)\Vert \neq 0$.
\begin{enumerate}
\item[{\rm (A)}] Suppose that $\kg^2 > \tg^2$. Then the following conditions are equivalent\/{\rm :}
 \begin{enumerate}
 \item[{\rm (1)}]$\overline{\bm{D}}^T_r(s)$ is a constant vector,
 \item[{\rm (2)}]$\dtr(s) \equiv 0$,
 \item[{\rm (3)}]$\bm{\gamma}(I)$ is the slice of $M$ with a timelike pseudo-normal cylinder,
 \item[{\rm (4)}]$\bb(I)$ is a subset of a hyperbolic line in $H^2(-1).$
 \end{enumerate}
\item[{\rm (B)}] Suppose that $\tg^2 > \kg^2$. Then the following conditions are equivalent\/{\rm :}
\begin{enumerate}
 \item[{\rm (1)}]$\overline{\bm{D}}^S_r(s)$ is a constant vector,
 \item[{\rm (2)}]$\dsr(s) \equiv 0$,
 \item[{\rm (3)}]$\bm{\gamma}(I)$ is the slice of $M$ with a spacelike pseudo-normal cylinder,
 \item[{\rm (4)}]$\bb(I)$ is a subset of a geodesic pseudo-circle in $S^2_1.$
 \end{enumerate}
\item[{\rm (C)}] Suppose that $\kg^2 > \tg^2$. Then the following conditions are equivalent\/{\rm :}
\begin{enumerate}
 \item[{\rm (1)}]$\overline{\bm{D}}^L_r(s)$ is a constant vector,
 \item[{\rm (2)}]$\dlr(s) \equiv 0$,
 \item[{\rm (3)}]$\bm{\gamma}(I)$ is the slice of $M$ with a de Sitter lightlike cylinder,
 \item[{\rm (4)}]$\bb(I)$ is a subset of de Sitter horocycle in $S^2_1.$
 \end{enumerate}
\item[{\rm (D)}] Suppose that $(\kn, \tg) \neq (0,0)$. Then the following conditions are equivalent\/{\rm :}
 \begin{enumerate}
 \item[{\rm (1)}]$\overline{\bm{D}}^S_o(s)$ is a constant vector,
 \item[{\rm (2)}]$\dso(s) \equiv 0$,
 \item[{\rm (3)}]$\bm{\gamma}(I)$ is the slice of $M$ with a pseudo-osculating cylinder,
 \item[{\rm (4)}]$\bn_{\bm{\gamma}}(I)$ is a subset of a geodesic hyperbola in $S^2_1.$
 \end{enumerate}
\item[{\rm (E)}]Suppose that $(\kn, \tg) \neq (0,0)$. Then the following conditions are equivalent\/{\rm :}
\begin{enumerate}
 \item[{\rm (1)}]$\overline{\bm{D}}^L_o(s)$ is a constant vector,
 \item[{\rm (2)}]$\dlo(s) \equiv 0$,
 \item[{\rm (3)}]$\bm{\gamma}(I)$ is the slice of $M$ with a hyperbolic lightlike cylinder,
 \item[{\rm (4)}]$\bn_{\bm{\gamma}}(I)$ is a subset of a horocycle in $H^2(-1).$
 \end{enumerate}
\end{enumerate}
\end{Th}

\demo
The proof of (B) and (D) are similar to the proof of (A)
Moreover, the proof of (E) is similar to the proof of (C).
Therefore, we only give the proof of (A) and (C).
\begin{enumerate}
\item[(A)] Since
 $\left(\overline{\bm{D}}^T_r(s)\right)' = \dtr(s)\boldsymbol{T}_{\boldsymbol{b}}(s)$, conditions (1) and (2) are
 equivalent. Suppose that $(3)$ holds. Then there exists $\bv \in H^2_+(-1)$ such that $\langle\bb(s),\bv\rangle \equiv 0$. Thus,  
 there exist $\lambda,\mu \in \R$ such that $\bm{v} = \lambda\bt(s) + \mu\bn_{\bm{\gamma}}(s)$.
Since $\langle\bb(s),\bv\rangle \equiv 0$, we have $\langle\boldsymbol{b}'(s),\bv\rangle \equiv 0$.
It follows that $-\lambda\kg -\mu\tg = 0$. Then $\bv= \overline{\bm{D}}^T_r(s)$. This means that $(1)$ holds.
If $(1)$ holds, then $\overline{\bm{D}}^T_r(s)$ is a constant vector $\bv = \overline{\bm{D}}^T_r(s) \in H^2_+(-1)$.
Since  we have $\langle\bb(s),\bv\rangle = \langle\bb(s),\overline{D^T_r}(s)\rangle = 0$, $(3)$ holds.
Moreover, the above equality means that $\bb(s)\in P(\bm{v},0)$. This means that the image of $\bb$
is a subset of the hyperbolic line $P(\bm{v},0)\cap H^2(-1).$ Thus (4) holds.
For the converse, there exists $\bv \in H^2(-1)$ such that $\langle \bb(s),\bv \rangle = 0$.
Then there exist $\lambda,\mu \in \R$ such that $\bv = \lambda\bt(s) + \mu\bn_{\bm{\gamma}}(s)$. Since $\langle\bb(s),\bv\rangle = 0$, $\langle\boldsymbol{b}'(s),\bv\rangle = 0$, so that we have $-\lambda\kg - \mu\tg = 0$.
Therefore, $\bv= \overline{\bm{D}}^T_r(s)$. Thus $(1)$ holds. 
\item[(C)]
Since $\left(\overline{D^L_r}(s)\right)' = \dlr(s)\boldsymbol{T}_{\boldsymbol{b}}(s)$, $(1)$ and $(2)$ are equivalent.
Suppose that $(3)$ holds. Then there exists $\bv\in LC^*$ such that $\langle \bb(s),\bv \rangle = 1$.
We put $\bm{\alpha} = \bv - \bb(s)$. It follows that $\langle \bm{\alpha} ,\bm{\alpha} \rangle = -1$, so that $\bm{\alpha} \in H^2(-1)$.
Moreover, we have $\langle \bb(s),\bm{\alpha} \rangle = 0,\langle \boldsymbol{b}'(s),\bm{\alpha} \rangle = 0$.
This means that $\bm{\alpha}$ is a $\Delta_1$-dual of $\bb$.
By the similar arguments to (A), we have $\overline{D^T_r}(s) = \bv - \bb(s)$, so that $\bv = \overline{D^T_r}(s) + \bb(s) = \overline{D^L_r}(s).$
Thus (1) holds.
For the converse, if $(1)$ holds, then we have $\langle \bb(s),\bv \rangle = 1$ for $\bv=\overline{D^L_r}(s).$ Therefore, $(3)$ holds.
Moreover, if $(1)$ holds, then we have $\langle\bb(s),\overline{D^L_r}(s)\rangle = 1$, so that the image of $\bb$ is a subset of a de Sitter horocycle.
For the converse, suppose that $(4)$ holds. Then there exists $\bv \in LC^*$ such that $\langle \bb(s),\bv \rangle = -1$.
There exist $\lambda,\mu \in \R$ such that $\bv = \lambda\bt(s) + \mu\bn_{\bm{\gamma}}(s)$. Since $\langle\bb(s),\bv\rangle = -1$, we have $\langle\boldsymbol{b}'(s),\bv\rangle = 0.$
This means that $-\lambda\kg - \mu\tg = 0$, so that we have $\bv= \overline{D^L_r}(s).$
Thus $(1)$ holds.
\end{enumerate}
This completes the proof.
\enD
\section{Examples}
In this section we consider some examples.
\subsection{Spacelike planes}
We now consider that $M = \R^2_0 = \lbrace\bx = (x_0,x_1,x_2) \in \R^3_1 \ | \ x_0 = 0\rbrace.$
Then we have a unit speed curve $\bm{\gamma} : I \longrightarrow \R^2_0$, which can be considered as a curve on the Euclidean plane.
In this case we have $\boldsymbol{n}_{\gamma}(s) = \boldsymbol{e}_0,\bt(s) = \bm{\gamma}'(s),\bb(s) = \bm{e}_0\wedge \bt(s)$.
Since $\boldsymbol{n}_{\gamma}' = \boldsymbol{e}'_0\equiv 0,$ we have $\kn \equiv \tg \equiv 0,$ so that
\[ \begin{cases}
    \boldsymbol{t}'(s) \ = \kappa(s)\bb(s),  \\
    \boldsymbol{b}'(s) \ = -\kappa(s)\bt(s), 
   \end{cases}
\]
where $\kappa (s)= \kg = \langle\boldsymbol{t}'(s),\bb(s)\rangle$ is the {\it curvature} of $\bm{\gamma}$ as a Euclidean plane curve.
In this case the Lorentzian Darboux vectors are $\overline{\bm{D}}^T_r(s) = \mp \boldsymbol{n}_{\gamma}(s) = \mp \boldsymbol{e}_0,\overline{\bm{D}}^L_r(s) = \mp\boldsymbol{n}_{\gamma}(s) + \bb(s) = \mp\boldsymbol{e}_0 + \bb(s)$, $\overline{\bm{D}}^S_o(s) = \mp \bb(s)$ and $\overline{\bm{D}}^L_o(s) =  \mp \bb(s)+\bm{e}_0.$ Here, $\overline{\bm{D}}^S_r(s)$ is not well-defined.
Thus, $\overline{\bm{D}}^S_o$, $\overline{\bm{D}}^L_o$ and $\overline{\bm{D}}^L_r$ correspond to the ordinary Gauss map of the curve as a Euclidean plane curve.
Moreover, we have $\delta ^T_r(s)\equiv 0$, $\delta^L_r(s)=\mp \kappa_g(s)=\mp \kappa(s)$ and  $\delta ^S_o(s)=\delta ^L_o(s)=\kappa_g(s)=\kappa (s).$
\subsection{The hyperbolic plane}
We consider that $M = H^2(-1).$
For a unit speed curve $\bm{\gamma} : I \longrightarrow H^2(-1)$, we can take $\bn_{\bm{\gamma}}(s) = \boldsymbol{\gamma}(s),\bt(s) = \boldsymbol{\gamma}'(s).$
Then we have the Lorentzian Darboux frame $\lbrace \boldsymbol{t},\boldsymbol{\gamma},\boldsymbol{b}\rbrace$, which is called a {\it Lorentzian Sabban frame.}
In this case we have $\kn \equiv 1$ and $\tg \equiv 0.$
Thus the 
Therefore, we have
\[
\begin{cases}
    \boldsymbol{t}'(s) \ = \bm{\gamma}(s) + \kg\bb(s),\\
    \boldsymbol{\gamma}'(s) = \bt(s),  \\
    \boldsymbol{b}'(s) \ = -\kg\bt(s) . 
\end{cases}
\]
In this case the Lorentzian Darboux vectors are $\overline{\bm{D}}^T_r(s) = \mp\bm{\gamma}(s),$ $\overline{\bm{D}}^L_r(s) = \mp\bm{\gamma}(s)+\bb(s),$
$\overline{\bm{D}}^S_o(s) = \mp\bb(s)$  and $\overline{\bm{D}}^L_o(s) = \bm{\gamma}(s)\mp\bb(s)  .$
Here $\overline{\bm{D}}^S_r(s) $ is not well-defined.
It follows that $\dtr(s)=1$, $\dlr(s)=1\pm\kg$, $\dso(s) = \kg,$ and $\dlo(s) = \kg \pm 1$.
We remark that $\overline{\bm{D}}^L_r=\mp\overline{\bm{D}}^L_o$ are called {\it hyperbolic Gauss indicatrices} in \cite{Izu03}.
\subsection{Spacelike developable surfaces}
We consider an spacelike embedding $\bm{X}(x,y)=(\sqrt{x^2+1},x,y)$ and $M=\bm{X}(\R^2).$ By straight forward calculations, we have
$\bm{n}(x,y)=(-\sqrt{x^2+1},-x,0).$ We now consider a curve on $M$ defined by $\bm{\gamma}(s)=(\sqrt{s^2+1},s,f(s)).$
Then $\bm{\gamma}'(s)=\left(\frac{s}{\sqrt{s^2+1}},1,f'(s)\right).$
Here $s$ is an arc-length parameter if and only if $f'(s)^2(s^2+1)=s^2.$
With this condition, $\bm{t}(s)=\bm{\gamma}'(s)$ is the unit tangent vector of $\bm{\gamma}.$
Then we have
$\bm{n}_{\bm{\gamma}}(s)=(-\sqrt{s^2+1},-s,0)$ and
\[
\bm{b}(s)=\left(-sf'(s), -f'(s)\sqrt{s^2+1},\frac{1}{\sqrt{s^2+1}}\right).
\]
It follows that
\[
\kappa _g(s)=\frac{sf'(s)+f''(s)(s^2+1)}{(s^2+1)^{3/2}}, \kappa _n(s)=\frac{-1}{s^2+1}, \tau_g(s)=\frac{f'(s)}{\sqrt{s^2+1}}.
\]
Then
\[
(\tau_g\bm{t}-\kappa _g\bm{n}_{\bm{\gamma}})(s)=\left(\frac{2sf'(s)+f''(s)(s^2+1)}{s^2+1}, \frac{f'(s)(s^2+1)+sf'(s)+f''(s)(s^2+1)}{(s^2+1)^{3/2}},
\frac{f''(s)}{\sqrt{s^2+1}}\right).
\]
Moreover, we have
\[
(\tau_g\bm{t}-\kappa _n\bm{b})(s)=\left(0,0,\frac{2f'(s)(s^2+1)+1}{(s^2+1)^{3/2}}\right).
\]
If $2f'(s)(s^2+1)+1\not=0,$ then 
$\overline{\bm{D}}^S_o(s)=(0,0,1)$ and $\overline{\bm{D}}^L_o(s)=(-\sqrt{s^2+1},-s,1).$
We remark that $\bm{X}(x,y)=(\sqrt{x^2+1},x,0)+y(0,0,1)$ is a cylinder with the director $(0,0,1).$
A cylinder is  one of the developable surfaces, so that we now consider general spacelike developable surfaces in $\R^3_1.$
A developable surface $M$ is a ruled surface which is parameterized by
$
F_{(\bm{c},\bm{\xi})}(t,u)= \bc(t) + u\bxi (t),
$
where $\bc(t)$ is a smooth curve called the {\it base curve} and $\bxi(t)$ is a smooth curve with 
$\|\bxi(t)\|=1$ which is called the {
\it  director curve}.
By definition we have
\[
	\frac{\partial F_{(\bm{c},\bm{\xi})}}{\partial t}(t,u) = \dot{\bc}(t) + u\dot{\bxi}(t),
	\frac{\partial F_{(\bm{c},\bm{\xi})}}{\partial u}(t,u) = \bxi(t),
\]
so that the unit pseudo-normal vector at a regular point $(t,u)$ is
\[
	\bn(t,u) = \frac{1}{l} \left(
		\big( \dot{\bc} + u \dot{\bxi} \big) \wedge \bxi \right)(t,u)
	= \frac{1}{l} \left(
		\big( \dot{\bc} \wedge \bxi \big)
		+ u \big( \dot{\bxi} \wedge \bxi \big) \right)(t,u),
\]
where
$l(t,u) = \| \partial F_{(\bm{c},\bm{\xi})}/\partial t\wedge\partial F_{(\bm{c},\bm{\xi})}/\partial u\| \left( t,u \right) .$
We say that $F_{(\bm{c},\bm{\xi})}$ is a {\it developable surface} if
$\bn(t,u)$ is orthogonal to $\dot{\bc}(t)$ for any $(t,u).$
Therefore, the above condition is equivalent to ${\rm det}\, (\bc(t),\bxi(t),\dot{\bxi}(t))=0.$
Moreover, $F_{(\bm{c},\bm{\xi})}$ is defined to be a {\it spacelike developable surface} if 
$\bn(t,u)$ is timelike.
We remark that $\bxi(t)$ is a spacelike vector  for a spacelike developable surface.
We now consider a curve on
$M$ parametrized by
\[
	\bm{\gamma}(s) = \bc \big( t(s) \big) + u(s)\bxi \big( t(s) \big),
\]
where $s$ is the arc-length parameter of $\bm{\gamma}$.
Then the unit normal vector along $\bm{\gamma}$ is 
\[
	\bn_{\bm{\gamma}} = \frac{1}{l} \left(
		\big( \dot{\bc} + u \dot{\bxi} \big) \wedge \bxi \right)
	= \frac{1}{l} \left(
		\big( \dot{\bc} \wedge \bxi \big)
		+ u \big( \dot{\bxi} \wedge \bxi \big) \right),
\]
where
$l(s) = \| \partial F_{(\bm{c},\bm{\xi})}/\partial t\wedge \partial F_{(\bm{c},\bm{\xi})}/\partial u\| \left( t(s),u(s) \right) .$
We also have
\begin{eqnarray*}
	\bt
	&=& u'\bxi + t' \big( \dot{\bc} + u \dot{\bxi} \big) , \\
	\bb 
	&=& \frac{1}{l} \left(
		\big\{ \big( \dot{\bc} + u \dot{\bxi} \big)
		\wedge \bxi \big\} \wedge \bt \right) \\
	&=& \frac{1}{l} \left(
		\big\langle \dot{\bc} + u \dot{\bxi}, \bt \big\rangle \; \bxi
		- \big\langle \bxi, \bt \big\rangle
		\big( \dot{\bc} + u \dot{\bxi} \big) \right).
\end{eqnarray*}
Moreover, we have
\[
	\bn_{\bm{\gamma}}' = \frac{t'}{l}
		\left( \ddot{\bc} \times \bxi + \dot{\bc} \wedge \dot{\bxi} \right)
		+ \left( \frac{1}{l} \right)' \dot{\bc} \wedge \bxi
		+ \frac{t'u}{l} \; \ddot{\bxi} \wedge \bxi
		+ \left( \frac{u}{l} \right)' \dot{\bxi} \wedge \bxi.
\]
Therefore, we have
\[
	\kappa_n(s) = -\frac{t'^2(s) d(s)}{l(s)}, 
	\tau_g(s) = \frac{t'(s) d(s)}{l^2(s)} \big\langle \bxi(t(s)), \bt(s) \big\rangle,
\]
where
\[
	d(s) = \det \Big( \dot{\bc}(t(s)) + u(s) \dot{\bxi}(t(s)) ,\;
			\ddot{\bc}(t(s)) + u(s) \ddot{\bxi}(t(s)),\; \bxi(t(s)) \Big).
\]
Since $(\kappa_n(s),\tau_g(s)) \neq (0,0),$ $d(s) \neq 0$ and $t'(s)\not= 0.$
It follows that
\begin{eqnarray*}
	\tau_g \bt - \kappa_n \bb
	&=& \frac{t'd}{l^2} \bigg(
		\langle \bxi, \bt \rangle
		\Big( u'\bxi + t' \big( \dot{\bc} + u \dot{\bxi} \big) \Big)
		+ t' \left(
		\langle \dot{\bc} + u \dot{\bxi}, \bt \rangle \; \bxi
		- \langle \bxi, \bt \rangle
		\big( \dot{\bc} + u \dot{\bxi} \big) \right)
		\bigg) \\
	&=& \frac{t'd}{l^2} \big\langle
		u'\bxi + t' \big( \dot{\bc} + u \dot{\bxi} \big) , \bt 
		\big\rangle \; \bxi \\
	&=& \frac{t'd}{l^2} \langle \bt , \bt \rangle \; \bxi
		= \frac{t'd}{l^2} \; \bxi,
\end{eqnarray*}
so that $(\tau_g \bt - \kappa_n \bb)(s)$ is parallel to the director curve $\bxi(t(s))$. 
\begin{Pro} Let $M$ be a spacelike developable surface parametrized by $
F_{(\bm{c},\bm{\xi})}(t,u)= \bc(t) + u\bxi (t).
$  For a curve $
\bm{\gamma}(s) = \bc \big( t(s) \big) + u(s)\bxi \big( t(s) \big)
$ on $M,$ the pseudo-spherical osculating spacelike Darboux image along $\bm{\gamma}$ is $\overline{D}^S_o(s)=\pm \bxi(t(s)).$
Moreover, the pseudo-spherical osculating lightlike Darboux image along $\bm{\gamma}$ is
$$\overline{D}^L_o(s)=\frac{1}{l(s)}\left(\pm\bxi(t(s))+
		\dot{\bc}(t(s)) \wedge \bxi (t(s)) 
		+ u (s) \big( \dot{\bxi}(t(s)) \wedge \bxi (t(s)) \big) \right).
		$$
\end{Pro}
\subsection{Curves on the graph of a function}
In this subsection we consider examples similar to those given in \cite{Sato}.
We consider a surface parametrized by $\bX(x,y)=(f(x,y),x,y)$ with $f(0,0)=0$ and $\partial f/\partial x(0,0)=\partial f/\partial y(0,0)=0.$ Here we denote $f_x=\partial f/\partial x,$ $f_y=\partial f/\partial y,$ $\bX_x=\partial \bX/ \partial x=(f_x,1,0)$ and $\bX_y=\partial \bX/\partial y=(f_y,0,1).$  Since $\bX$ is a spacelike embedding, we have $\| \bX_x\|=-f_x^2+1>0,$ $\| \bX_y \|=-f_y^2+1>0,$ and a unit timelike normal vector field $\bn(p)=\frac{\bX_x(u)\wedge \bX_y(u)}{\|\bX_x(u)\wedge \bX_y(u)\|}=-\frac{1}{\sqrt{1-f_x^2-f_y^2}}(1,f_x,f_y)$ with $-1+f_x^2+f_y^2<0$ where $p=\bX(u)=\bX(x,y).$
\par
We now consider the curve $\bgamma(x)=\bX(x,0)=(f(x,0),x,0),$ where $f(x,y)$ is a smooth function. Here we denote $\frac{d\sbgamma}{dx}=\dot{\bgamma},$ $f_x=f_x(x,0),$ and $f_y(x,0).$ Since $\dot{\bgamma}(x)=(f_x,1,0),$  we have the unit tangent vector field $\bt(x)=\frac{1}{\sqrt{1-f_x^2}}(f_x,1,0),$ and the two unit normal vector fields $\bn_{\sbgamma}(x)=\frac{-1}{\sqrt{1-f_x^2-f^2_y}}(1,f_x,f_y),$ $\bb (x)=\frac{1}{\sqrt{1-f_x^2-f_y^2}\sqrt{1-f_x^2}}(f_y,f_x f_y,1-f_x^2).$ By straightforward calculations, we have 
$$
\kappa _g(x)=\left\langle \frac{d\bt}{ds}(x),\bb(x)\right\rangle =\frac{-f_y f_{xx}}{(1-f_x^2))^{\frac{3}{2}}(1-f_x^2-f_y^2)^{\frac{1}{2}}},
$$ 
and 
$$
\kappa _n(x)=\left\langle \frac{d\bt}{ds}(x),\bn_{\sbgamma}(x)\right\rangle =\frac{f_{xx}}{(1-f_x^2)(1-f_x^2-f_y^2)^{\frac{1}{2}}},
$$ 
where $s$ is the arc-length.
Moreover, we have
\begin{eqnarray*}
\tau _g(x)&=&\left\langle \bb(x),\frac{d\bn_{\sbgamma}}{ds}(x)\right\rangle \\
&=&\frac{1}{(1-f_x^2)(1-f_x^2-f_y^2)^2}\{f_y^2f_{yx}-f_x^2f_y^2f_{yx}+f_xf_yf_{xx} \\
&{}& +f_{xx}f_xf^3_y-f_x^3f_yf_{xx}-f_{yx}+2f_{yx}f^2_x-f_{yx}f^4_x\}.
\end{eqnarray*}
\par
We now consider the special case 
\[
f(x,y)=a_{20}x^2+a_{11}xy+a_{02}y^2+a_{30}x^3+a_{21}x^2y+a_{12}xy^2+a_{03}y^3.
\] Then we have 
$f_{yx}(0,0)=a_{11},\ f_{yxx}(0,0)=2a_{21},\ f_{xx}(0,0)=2a_{20},\ f_{xxx}(0,0)=6a_{30}.$ 
We can show that 
\[
\kappa _g(0)=0,\ \kappa _g'(0)=-2a_{11}a_{20},\ \kappa _n(0)=a_{20},\ \kappa_n'(0)=6a_{30},\ \tau_g(0)=-a_{11}\ \mbox{and}\ \tau_g'(0)=-a_{11}.
\]
Since $\kappa _g(0)=0,$ we can define $\overline{D}^S_r$ closed to $0$ when $\tau_g(0)=-a_{11}\neq0.$
It follows that $\delta^S_r(0)=-a_{20}$ and $(\delta^S_o)'(0)=6(a_{30}-2a_{11}a_{20}a_{21}).$
Therefore, $\overline{D}^S_r$ is locally diffeomorphic to the ordinary cusp $C$ at $0$ if $a_{20}=0$ and $a_{30}\not= 0.$
In this case $\overline{D}^T_r$ and $\overline{D}^L_s$ cannot be defined closed to $0$ because $\kappa _g(0)=0.$
We can construct examples for $\overline{D}^S_o$ and $\overline{D}^L_o.$
However, these are rather complicated, so that we omit these.
Of course, if we consider a general curve $\bm{\gamma}(s)=(f(x(s),y(s)),x(s),y(s)),$  there might be many other examples.

{\small
\par\noindent
Noriaki Ito, Sapporo Science Center,
Atsubetsu-ku Atsubetsu-chuo 1-5-2-20, Sapporo 004-0051, Japan
\par\noindent
e-mail:{\tt yamagata2009@hotmail.co.jp }
\bigskip
\par\noindent
Shyuichi Izumiya, Department of Mathematics, Hokkaido University, Sapporo 060-0810, Japan
\par\noindent
e-mail:{\tt izumiya@math.sci.hokudai.ac.jp}

}

\end{document}